\documentclass[12pt]{article}

\usepackage{amsmath,amssymb,amsfonts}
\usepackage{graphicx}
\usepackage{float}
\usepackage{booktabs}
\usepackage{siunitx}
\usepackage[numbers,sort&compress]{natbib}
\usepackage{hyperref}

\hypersetup{
    colorlinks=true,
    linkcolor=blue,
    urlcolor=blue,
    citecolor=blue
}

\title{Numerical Behavior of the Riemann Zeta Function Using Real-to-Complex Conversion Without Gram Points or Bracketing}

\author{Jacob Orellana Real\\
Independent Researcher, Ecuador\\
\texttt{jacoboreore@gmail.com}
}

\date{December 2025}

\begin{document}
\maketitle

\begin{center}
\textbf{Preprint submitted to arXiv}\\
This version: December 2025
\end{center}
\vspace{1em}

\begin{abstract}
The Riemann zeta function $\zeta(s)$ and its zeros lie at the core of modern analytic number theory. 
On the critical line $\Re(s) = \tfrac12$, the Hardy function
\[
  Z(t) = e^{i\theta(t)}\,\zeta\!\left(\tfrac12 + it\right)
\]
is real valued and encodes the nontrivial zeros as sign changes. 
Classical large-scale computations rely heavily on Gram points, Gram intervals, and bracketing strategies for locating zeros.

In this paper we develop a different approach for locating zeros of $Z(t)$ based on a real-to-complex conversion
\[
  t = \sqrt{N - \tfrac14},
\]
which parametrizes large heights $t$ by a real variable $N$. 
We combine this parametrization with the Riemann--Siegel formula and Gabcke-type remainder estimates, and introduce a 
``valley scanner'' that detects local minima of $\lvert Z(t)\rvert$ along the resulting trajectory in $t$-space
rather than traditional sign-change accounting, and offers a flexible computational tool for further large-scale experiments 
on the critical line.
\end{abstract}

\bigskip


\medskip


\section{Introduction}

The numerical exploration of the Riemann zeta function $\zeta(s)$ on the
critical line $s = \tfrac{1}{2} + it$ remains a central component of analytic
number theory. Large-scale computations of its nontrivial zeros have
historically relied on Gram points, Gram intervals, and associated bracketing
techniques developed by Titchmarsh, van de Lune and collaborators, and,
more recently, Odlyzko. These methods have proven effective for locating and
verifying zeros at extremely large heights, but they also introduce structural
constraints: the distribution of Gram points becomes increasingly irregular,
Gram intervals may fail to bracket a zero, and the analysis often requires
delicate case-by-case logic.

This paper introduces an alternative numerical strategy based on a
real-to-complex transformation that maps a real parameter $N$ to a
corresponding height $t$ on the critical line. Although the transformation is
elementary, it induces a geometric structure that interacts favorably with the
oscillatory behavior of the Hardy function $Z(t)$. In particular, plotting
$\lvert Z(t)\rvert$ under this parametrization reveals a landscape of mountains
and valleys whose local minima naturally correspond to neighborhoods of zeta
zeros. This observation motivates the construction of a \emph{valley scanner}:
a procedure that detects candidate minima of $\lvert Z(t)\rvert$ and uses these
points as starting positions for Newton refinements.

The computational framework developed here combines:
\begin{itemize}
    \item the real-to-complex parametrization $N \mapsto t$,
    \item high-precision evaluation of the Riemann--Siegel formula,
    \item Gabcke-type bounds for the remainder term,
    \item a safeguarded Newton refinement.
\end{itemize}

Extensive computations were performed for heights up to
$t \approx 10^{20}$, with results compared against existing zero datasets. The
method consistently reproduces known zeros and displays stable behavior well
beyond the range of public tables. Because it does not rely on Gram intervals,
the approach offers a geometry-driven alternative to classical bracketing and
may provide a complementary viewpoint for future large-scale studies of the
zeta function.

The algorithms, datasets, and validation tools used in this study have been
publicly released for reproducibility through a persistent DOI. All numerical
results reported here were generated using these open resources.


\section{Mathematical Background}
\label{sec:background}

\subsection{The Riemann Zeta Function}

The Riemann zeta function is defined for complex arguments
$s = \sigma + it$ with $\sigma > 1$ by the absolutely convergent series
\[
\zeta(s) = \sum_{n=1}^{\infty} \frac{1}{n^s}.
\]
Through analytic continuation, this definition extends to all
$s \neq 1$. A central identity governing its structure is the
functional equation
\[
\zeta(s)
=
2^s \pi^{\,s-1}
\sin\!\left(\frac{\pi s}{2}\right)
\Gamma(1-s)\,\zeta(1-s),
\]
which reflects values across the critical line. The nontrivial zeros
of $\zeta(s)$ lie in the strip $0 < \sigma < 1$.
The Riemann Hypothesis asserts that all such zeros satisfy
\[
\Re(s) = \tfrac{1}{2}.
\]

\subsection{Relation to Prime Numbers}

The zeta function encodes the multiplicative structure of the positive
integers through Euler's product
\[
\zeta(s)
= 
\prod_{p \,\text{prime}}
\frac{1}{1 - p^{-s}}, 
\qquad
\Re(s) > 1.
\]
Its zeros influence the distribution of primes via the explicit formula
relating $\pi(x)$ to the nontrivial zeros
$\rho = \tfrac{1}{2} + i t_{\rho}$:
\[
\pi(x)
\sim
\operatorname{Li}(x)
-
\sum_{\rho} \operatorname{Li}(x^{\rho})
+
\text{(lower-order terms)}.
\]
This deep connection motivates the search for efficient and accurate methods
for computing zeros at large heights.

\subsection{The Hardy Function $Z(t)$}

On the critical line, the real-valued Hardy function is defined by
\[
Z(t)
=
e^{i\theta(t)} \,
\zeta\!\left(\tfrac{1}{2} + it\right),
\]
where the Riemann--Siegel theta function is
\[
\theta(t)
=
\arg\!\left(
\Gamma\!\left(\tfrac{1}{4} + \tfrac{it}{2}\right)
\right)
-
\frac{t}{2} \log \pi.
\]
The factor $e^{i\theta(t)}$ ensures that $Z(t)$ is real for real $t$.
The magnitude $\lvert Z(t)\rvert$ exhibits a characteristic
mountain--valley structure, and its zeros correspond precisely to the
nontrivial zeros of $\zeta(s)$ on the critical line.

\subsection{Average Zero Spacing and Density}

The Riemann--von Mangoldt formula states that the number of nontrivial
zeros with imaginary part in $(0, T)$ is
\[
N(T)
=
\frac{T}{2\pi}
\log\!\left(\frac{T}{2\pi}\right)
-
\frac{T}{2\pi}
+
O(\log T).
\]
Thus, the average spacing between consecutive zeros near height $T$ is
\[
\Delta T_{\text{avg}}
\approx
\frac{2\pi}{\log\!\left(T / 2\pi\right)}.
\]
This estimate provides a meaningful benchmark for assessing the density
and distribution of numerically computed zeros.

\subsection{Known Computational Approaches}
\label{subsec:known-approaches}

The numerical evaluation of $Z(t)$ and the computation of its zeros have led
to a wide range of techniques, from asymptotic expansions to FFT-based
algorithms.

\subsubsection*{The Riemann--Siegel Formula}

For large $t$, the Dirichlet series representation of $\zeta(s)$ converges
too slowly to be practical. The Riemann--Siegel formula provides the asymptotic
expansion
\[
Z(t)
=
2 \sum_{n=1}^{N}
\frac{\cos(\theta(t) - t\log n)}{\sqrt{n}}
+
R(t),
\qquad
N = \left\lfloor \sqrt{\frac{t}{2\pi}} \right\rfloor,
\]
with a small remainder term $R(t)$. Gabcke's refinements improve the evaluation
of the remainder and help ensure numerical stability at large heights.

\subsubsection*{Gram Points}

Gram points are values $g_n$ defined implicitly by
\[
\theta(g_n) = n\pi.
\]
Empirically, zeros of $Z(t)$ tend to alternate with Gram points, although
Gram's law fails infinitely often. Nonetheless, Gram points remain a classical
tool for locating sign changes.

\subsubsection*{Turing's Method}

Turing introduced a verification procedure \cite{Turing} that ensures no zeros
are missed in an interval by comparing the number of sign changes with the
theoretical count predicted by the Riemann--von Mangoldt formula. This remains
the standard completeness check in computational work.

\subsubsection*{Odlyzko--Sch\"onhage Algorithm}

At very large heights, the Odlyzko--Sch\"onhage method \cite{OdlyzkoSchonhage}
uses FFT-based acceleration to evaluate $Z(t)$ at many points simultaneously.
This approach enabled computations of millions of zeros far beyond $10^{20}$
and remains the state of the art for ultra-high-range verification.

\subsubsection*{Limitations at High $t$}

As $t$ increases:
\begin{itemize}
    \item the oscillation of $Z(t)$ intensifies, increasing the required working precision;
    \item the Riemann--Siegel expansion becomes increasingly sensitive to remainder terms;
    \item FFT-based methods face memory-scaling limitations and difficulties integrating arbitrary precision.
\end{itemize}

These constraints motivate alternative formulations aimed at numerical
stability and conceptual simplicity, such as the real-to-complex conversion
and valley-scanning methods developed in this work.

\section{Proposed Numerical Method}
\label{sec:method}

\subsection{Motivation and Research Shift}

Although the present work focuses on the Riemann zeta function and the
numerical detection of its zeros, the underlying method originated from
an independent study of semiprime numbers. The original goal was to
analyze algebraic relationships between a semiprime $N = p_1 p_2$ and its
constituent prime factors by examining symmetry structures embedded in
associated quadratic forms.

Given
\[
N = p_1 p_2, \qquad S = p_1 + p_2,
\]
the classical quadratic equation
\[
x^2 - Sx + N = 0
\]
has roots $p_1$ and $p_2$, with solutions
\[
x = \frac{S \pm \sqrt{S^2 - 4N}}{2}.
\]
Although algebraically natural, this framework did not lead to an
efficient computational method for factor analysis and was eventually
set aside. Nonetheless, the quadratic structure suggested a viewpoint
that proved unexpectedly fruitful in the context of zeta-function
computations.

\subsection{N-to-Complex Conversion}

Rather than focusing directly on the unknown primes, the analysis shifted
toward their arithmetic mean:
\[
m = \frac{p_1 + p_2}{2}.
\]
Since $p_2 = 2m - p_1$, substituting into $N = p_1 p_2$ yields
\[
N = p_1(2m - p_1),
\]
which rearranges to the quadratic equation
\[
p_1^2 - 2mp_1 + N = 0.
\]
The quadratic formula gives
\[
p_1 = m \pm \sqrt{m^2 - N}.
\]

Originally the analysis focused on real-valued geometric behavior such as
the area under the quadratic curve, but no stable correlation with the
semiprime structure emerged. This led to an alternative interpretation:
whenever $m^2 < N$, the expression naturally produces complex solutions,
\[
p_1 = m \pm i\sqrt{N - m^2},
\]
suggesting that the quadratic form defines a mapping from each real $N$
to a complex value.

A natural question is: \emph{What happens if the mean is fixed at a
specific value?} Setting
\[
m = \tfrac{1}{2}
\]
aligns the construction with the critical line $\Re(s) = \tfrac{1}{2}$ of
the Riemann zeta function. With this choice one obtains
\[
p(N) = \frac{1}{2} \pm i\sqrt{N - \frac{1}{4}},
\]
leading to the real-to-complex mapping
\[
N \longmapsto s(N)
=
\frac{1}{2} \pm i\sqrt{N - \tfrac{1}{4}}.
\]
Each real input $N$ maps to a point on the critical line with imaginary
part determined by the square-root term. Numerical exploration of this
mapping revealed structured behavior in the corresponding evaluations of
the Hardy $Z$-function, inspiring the development of the valley-scanner
algorithm.

\subsection{Valley Scanner Algorithm}

The term \emph{valley scanner} refers to the numerical procedure used to
identify zeros of the Hardy function by detecting local minima of
$\lvert Z(t)\rvert$ along a parametrized traversal of the critical line.
Starting with
\[
s = m \pm i\sqrt{N - m^2}, \qquad m = \tfrac{1}{2},
\]
we have
\[
t = \Im(s) = \sqrt{N - m^2}, \qquad
N = t^2 + m^2.
\]
Thus increments in $N$ induce increments in $t$, providing a simple
mechanism for scanning the values of $\lvert Z(t)\rvert$ across a
continuous range.

The objective is to evaluate $\lvert Z(t)\rvert$ for successive values of
$N$ and inspect the resulting numerical landscape.

An initial Python implementation was used for exploratory experiments.
The formatted version of the script is available at:\newline

\url{https://github.com/zjore/z-research/blob/main/playground.py}.

\begin{figure}[H]
\centering
\includegraphics[width=0.85\textwidth]{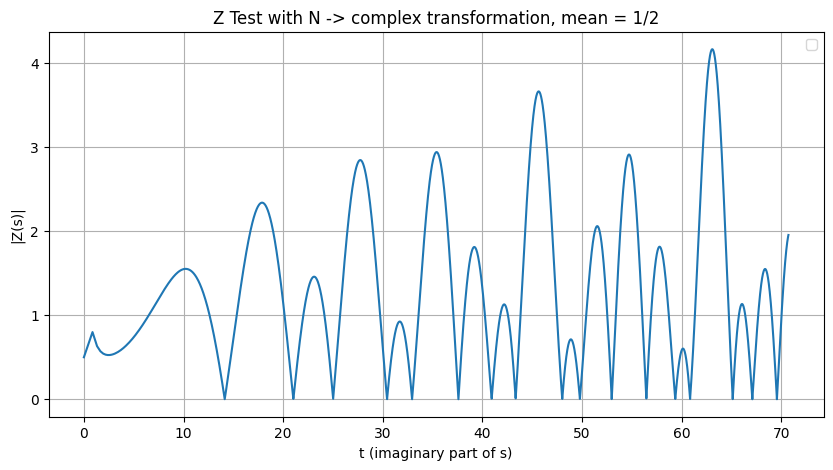}
\caption{Visualization of valley detection through local minima of
$\lvert Z(t)\rvert$.}
\label{fig:valley}
\end{figure}

A comparison with the first nontrivial zeros of the Riemann zeta function is
shown in Table~\ref{tab:riemann-zeros}.

\begin{table}[H]
\centering
\caption{First nontrivial zeros of the Riemann zeta function.}
\begin{tabular}{ll}
\toprule
\textbf{Index} & \textbf{Imaginary part $t_n$}  \\
\midrule
1  & 14.1347251417346937904572519835625 \\
2  & 21.0220396387715549926284795938969 \\
3  & 25.0108575801456887632137909925628 \\
4  & 30.4248761258595132103118975305840 \\
5  & 32.9350615877391896906623689640747 \\
6  & 37.5861781588256712572177634807053 \\
7  & 40.9187190121474951873981269146334 \\
8  & 43.3270732809149995194961221654068 \\
9  & 48.0051508811671597279424727494277 \\
10 & 49.7738324776723021819167846785638 \\
11 & 52.9703214777144606441472966088808 \\
12 & 56.4462476970633948043677594767060 \\
13 & 59.3470440026023530796536486749922 \\
14 & 60.8317785246098098442599018245241 \\
15 & 65.1125440480816066608750542531836 \\
16 & 67.0798105294941737144788288965221 \\
17 & 69.5464017111739792529268575265547 \\
\bottomrule
\end{tabular}
\label{tab:riemann-zeros}
\end{table}

Visually, the resulting plot is self-explanatory. When the Riemann zeta
function $\zeta(s)$ is evaluated using the parametrization
$s = m \pm i\sqrt{N - m^{2}}$ with $m = \tfrac{1}{2}$, and with $N$
ranging over an arbitrary numerical interval, the computed values
consistently approach zero whenever the corresponding height $t$
is near a true zero of the zeta function. This behavior is fully
consistent with the theoretical expectation that the transformation
traces the critical line and that the magnitude of $\zeta(s)$ should
decrease in the vicinity of its nontrivial zeros.

Refining the range of $N$ by introducing decimal-level resolutions
produces progressively closer approximations to the actual zeros.
The accompanying Python script reproduces this behavior across any
chosen interval of $N$, though the computations become increasingly
demanding as the height $t$ grows.

The dataset obtained from this exploratory script is summarized below
and compared against confirmed zeros of the Riemann zeta function.

\begin{table}[H]
\centering
\caption{Comparison between computed and reference zeros of the Riemann zeta function.}
\begin{tabular}{lll}
\toprule
\textbf{$N$} & \textbf{Computed $t$} & \textbf{Reference $t$} \\
\midrule
200  & 14.13329402510257  & 14.13472514 \\
442  & 21.017849556983705 & 21.02203964 \\
626  & 25.014995502697975 & 25.01085758 \\
926  & 30.426140077242792 & 30.42487613 \\
1085 & 32.93554311074891  & 32.93506159 \\
1413 & 37.586566749305526 & 37.58617816 \\
1675 & 40.92370951        & 40.91871901 \\
1877 & 43.32147273581544  & 43.32707328 \\
2305 & 48.00781186432057  & 48.00515088 \\
2478 & 49.77700272214067  & 49.77383248 \\
2806 & 52.96933074902875  & 52.97032148 \\
3186 & 56.442448564887755 & 56.44624770 \\
3522 & 59.34433418617145  & 59.34704400 \\
3701 & 60.83378995        & 60.83177853 \\
4240 & 65.11336268386084  & 65.11254405 \\
4500 & 67.08017591        & 67.07981053 \\
4837 & 69.54674686856315  & 69.54640171 \\
\bottomrule
\end{tabular}
\label{tab:computed_vs_public}
\end{table}

\subsection{Identifying Zeros via Consecutive Minima of $\lvert Z(t)\rvert$}

Although Python provides a flexible environment for preliminary
experiments, the valley-scanning approach becomes computationally
intensive as the height $t$ increases. To explore larger regions of the
critical line, optimized implementations were developed in C++ with
arbitrary-precision arithmetic, together with parallel execution on AWS
EC2 instances.

The underlying idea remains unchanged: the magnitude $\lvert Z(t)\rvert$
exhibits an alternating pattern of peaks and valleys, and the nontrivial
zeros occur precisely at the bottoms of these valleys. The algorithm
therefore scans successive values of $N$ (and the corresponding $t$)
to detect transitions from descent to ascent in the numerical landscape.

\subsubsection*{Mountain walking and valley detection}

The valley-scanner algorithm follows a simple geometric principle:

\begin{itemize}
    \item Choose an initial height $t_0$.
    \item Convert $t_0$ to its corresponding value $N_0$ via $N = t^{2} + m^{2}$.
    \item Increment $N$ by one unit (or by a chosen step size).
    \item Compute the resulting $t_1$ and evaluate $\lvert Z(t_1)\rvert$.
    \item If $\lvert Z(t_1)\rvert < \lvert Z(t_0)\rvert$, the scan is descending; otherwise it is ascending.
    \item A change from descent to ascent indicates a local minimum and therefore a candidate zero.
    \item Each candidate is refined by a Newton root-finding method.
\end{itemize}

This mechanism serves as the foundation for the large-scale computations
described later in the paper.

\subsubsection*{Evaluation of the Hardy function}

The numerical evaluation of $Z(t)$ relies on several components:

\begin{itemize}
    \item The Riemann--Siegel formula
    \[
        Z(t)
        =
        2 \sum_{n=1}^{N}
        \frac{\cos(\theta(t) - t \log n)}{\sqrt{n}}
        + R(t),
    \]
    where $N = \lfloor \sqrt{t/(2\pi)} \rfloor$, and $R(t)$ is the
    remainder term estimated using Gabcke-type refinements.

    \item The phase function $\theta(t)$ is computed via a high-precision
    evaluation of $\log \Gamma(\tfrac{1}{4} + \tfrac{it}{2})$ using
    Stirling expansions:
    \[
        \theta(t)
        =
        \operatorname{Im}
        \left[
          \log \Gamma\!\left(\frac{1}{4} + \frac{it}{2}\right)
        \right]
        - \frac{t}{2}\log \pi.
    \]

    \item Zero candidates are refined using a safeguarded Newton method
    with a Pegasus-type fallback (a hybrid Newton--secant scheme).
    Computations are performed at 128--256 bits of precision to stabilize
    convergence in highly oscillatory regions.

    \item Validation is performed by comparing the refined zeros with
    independent datasets and by applying a local consistency test (the
    ``ball test'') across neighboring zeros.
\end{itemize}

\subsection{Parallelization and Cloud Execution}

The dominant computational cost in evaluating $Z(t)$ through the
Riemann--Siegel formula arises from the main summation
\[
    Z(t)
    =
    2\sum_{n=1}^{N}
      \frac{\cos(\theta(t) - t\log n)}{\sqrt{n}}
    + R(t).
\]
For fixed $t$, each term $\cos(\theta(t) - t\log n)/\sqrt{n}$ is
independent, making the summation highly amenable to parallelization.

We partition the summation interval into subranges
\[
S_i
=
\sum_{n = n_i}^{n_{i+1}-1}
  \frac{\cos(\theta(t) - t\log n)}{\sqrt{n}},
\]
and aggregate the partial sums. The remainder term $R(t)$ is evaluated
separately using Gabcke's expansion.

Only the main summation is parallelized; the correction term remains
serial. In practice this yields near-linear scaling with the number of
physical CPU cores. AWS EC2 instances are configured with OpenMP to
distribute work across all cores using deterministic reduction clauses,
ensuring numerical stability and reproducibility.

This parallel strategy remains effective even when
\( N = \sqrt{t/(2\pi)} \) reaches into the millions, as occurs for very
large heights.

\subsection{Why Spacing Heuristics (and Naive Prefiltering) Can Skip Zeros}
\label{sec:spacing-skip}

The Riemann--von Mangoldt formula
\begin{equation}
  N(T)
  =
  \frac{T}{2\pi}\log\!\Bigl(\frac{T}{2\pi}\Bigr)
  - \frac{T}{2\pi}
  + O(\log T)
  \label{eq:NT}
\end{equation}
implies an approximate zero density
\[
N'(T)
\approx
\frac{1}{2\pi}
\log\!\Bigl(\frac{T}{2\pi}\Bigr),
\]
and an average spacing
\begin{equation}
  s_{\mathrm{avg}}(T)
  \approx
  \frac{2\pi}{\log(T/(2\pi))}.
  \label{eq:avg-spacing}
\end{equation}

\paragraph{Key caveat.}
While \eqref{eq:avg-spacing} provides a reliable mean estimate, the
\emph{local} spacing between consecutive zeros can differ substantially.
Advancing $t$ by a deterministic step of size $s_{\mathrm{avg}}(T)$, or
rejecting points based on naive thresholds, may cause genuine zeros to
be skipped.

\paragraph{Concrete example (skipped zero).}
Two consecutive zeros verified independently satisfy
\[
t_1 = 30607946001.041439,
\qquad
t_2 = 30607946001.175073,
\]
yielding a gap
\[
\Delta_{\mathrm{true}}
=
t_2 - t_1
\approx
0.133634.
\]
The spacing estimate at \(T \approx t_1\) is
\[
s_{\mathrm{avg}}(t_1)
=
\frac{2\pi}{\log(t_1/(2\pi))}
\approx
0.281673,
\]
which overshoots the next zero by a wide margin.

\paragraph{Implication for the valley scanner.}
Early variants of the algorithm explored spacing-based stepping and
threshold-based prefilters. Both approaches resulted in missed zeros.
The final version therefore avoids prefilters and relies instead on
dense sampling with minima detection, followed by safeguarded
refinement.

\subsubsection*{Software and Computational Stack}

The computational workflow is built on:

\begin{itemize}
    \item \textbf{C++} with MPFR/MPC for high-precision evaluation of $Z(t)$.
    \item \textbf{AWS EC2} for large multi-core compute instances.
    \item \textbf{AWS Lambda} for batch orchestration and monitoring.
    \item \textbf{AWS S3} for persistent storage of logs and zero tables.
    \item \textbf{AWS Cognito} for secure authentication.
    \item \textbf{Docker} and AWS ECR for reproducible deployments.
\end{itemize}

An authenticated web interface allows users to run valley scans and
visualize zeros up to heights near $10^{16}$. While the full source code
is not yet archived, all numerical datasets used in this work are
publicly available via Zenodo and GitHub. A public Docker image is
provided for reproducibility; see the project documentation:\newline

\url{https://github.com/zjore/z-research/blob/main/README.md}


\section{Results}

This section presents representative outputs obtained by the
valley-scanner framework across several ranges of $t$, including large
heights where refinement is particularly challenging. The analysis of
average zero spacing and local fluctuations was cross-checked by
independent numerical scripts and symbolic tools to validate the
results.

\paragraph{Scope of datasets.}
The datasets shown here are intentionally concise and primarily
illustrative. The valley scanner is designed to detect zeros within any
specified range of $t$, and the batch-processing architecture supports
large-scale computation over much broader intervals, with feasibility
limited mainly by practical cloud-computing costs.

\begin{table}[H]
\centering
\small
\caption{Zero heights near $t \approx 1.1223\times 10^{9}$.}
\begin{tabular}{l}
\hline
\textbf{$t$ values} \\
\hline
1122334455.05585408179597553373792753809670000283010284579445795581 \\
1122334455.76031733677095186350835049571870227573769745536211492284 \\
1122334455.96549341079558436082594904200898020740502977379833207345 \\
1122334456.55122003455096567518805961663152718628678713782068430266 \\
1122334456.69117656282847361682640086578613693696115124531479812624 \\
1122334456.83867869196755030364329563494267434476218915992393583291 \\
1122334457.15371418683945191486723066707952951489968507756091862973 \\
1122334457.51390511989118717254611539418981231215406357632944032390 \\
1122334457.76838741510620113455239734916962707865127304851980893696 \\
1122334458.27155512243026569374168652041355212655972394951814858601 \\
1122334458.44139136841692632939134894613526508309356307629078128927 \\
1122334458.89501425212585279215913084722248689346999249199510723479 \\
1122334458.95374621822277785558114513758006054444953092856654909690 \\
1122334459.47180990048612714347328799643237077235245383445928401233 \\
1122334459.78629512841632863716858891532650291270784909335852103315 \\
1122334460.10222132816663933236459832453382800744621284674207091481 \\
1122334460.34922785530398779220541441366240746601424248348045312059 \\
1122334460.78924267475219010097101478269117727805242937830818188169 \\
1122334461.03129979223085870682939517578522137411650684496149763174 \\
1122334461.34105411839180789869245630301714029644280234966395722460 \\
1122334461.68097513340590425015663066588174538735332969362229273055 \\
1122334462.05653607553267581154859487598541969417453119764160088476 \\
1122334462.35999091678192331701334666903361418145334941591973792237 \\
1122334462.71045060502215878868517734546324575820479730409902219700 \\
1122334463.22570175761036362585961260605296584140311386148887315924 \\
\hline
\end{tabular}
\label{tab:t_values_1122334455}
\end{table}

\noindent
Full dataset: \\
\url{https://github.com/zjore/z-research/blob/main/datasets/refined_sample_1122334455.csv}

\begin{figure}[H]
    \centering
    \includegraphics[width=0.85\textwidth]{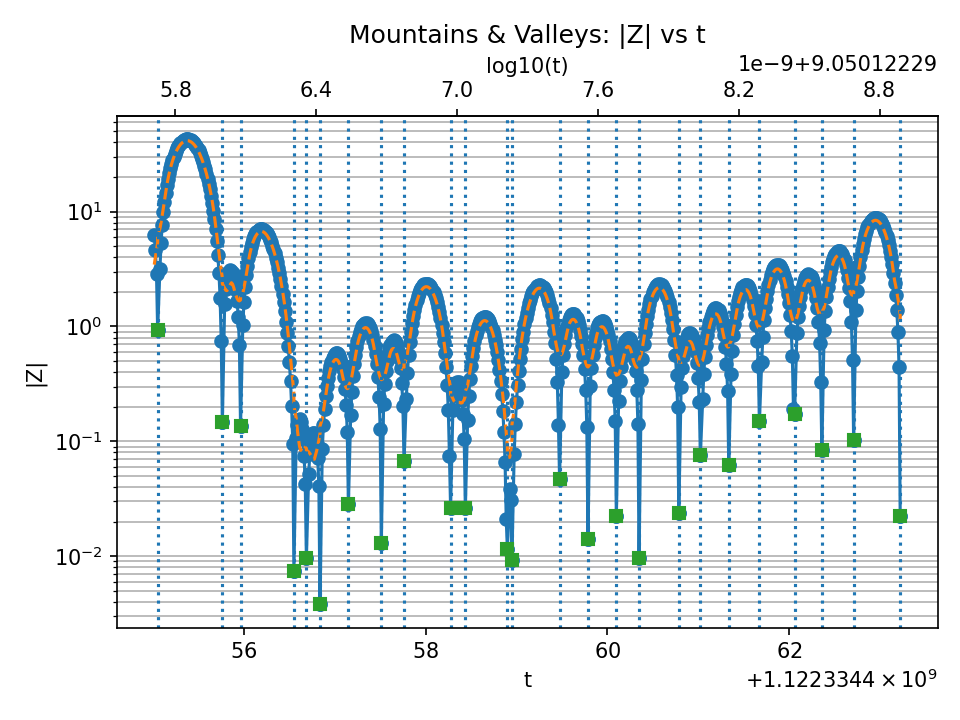}
    \caption{Valley-scanner output near $t \approx 1.1\times 10^9$; local minima of $\lvert Z(t)\rvert$ are marked prior to refinement.}
    \label{fig:valley_mid}
\end{figure}

\begin{table}[H]
\centering
\caption{Observed zero spacing vs.\ theoretical prediction at $t \approx 1.1223\times 10^{9}$.}
\begin{tabular}{ll}
\hline
\textbf{Parameter} & \textbf{Value / Description} \\
\hline
Mean height $T$ & $1.122334459\times 10^{9}$ \\
Theoretical formula
 & $\Delta_{\mathrm{avg}} \approx \dfrac{2\pi}{\log(T / 2\pi)}$ \\
$\log(T / 2\pi)$ & $19.84$ \\
Predicted $\Delta_{\mathrm{avg}}$ & $0.33068$ \\
Observed mean$(\Delta t)$ & $0.34041$ \\
Relative difference & $\approx 2.94\%$ \\
\hline
\end{tabular}
\label{tab:spacing_comparison_1122334455}
\end{table}

The observed mean spacing differs from the theoretical prediction by
about $2.9\%$, which is consistent with ordinary local fluctuations at
this height. This dataset contains twenty-four valid spacings (after
discarding the initial $\Delta t = 0$), providing a stable estimate of
the zero density.

\paragraph{Additional comments.}
\begin{itemize}
\item The narrowest local gap ($\Delta t \approx 0.05873$) reflects a
      brief clustering region.
\item The widest local gap ($\Delta t \approx 0.70446$) corresponds to
      a modest sparsity interval.
\item Across this range, the valley scanner maintained $\lvert Z(t)\rvert$
      residuals below $10^{-17}$ at refined zeros.
\end{itemize}

\begin{table}[H]
\centering
\small
\caption{Zero heights near $t \approx 7.7888\times 10^{13}$.}
\begin{tabular}{l}
\hline
\textbf{$t$ values} \\
\hline
77887788778877.0202277488171921884175366387426290737402889145896393 \\
77887788778877.2496796799624248241926909299697953447734769996091705 \\
77887788778877.3958394524380600874237302170777958741989520806110630 \\
77887788778877.4642294878030461471735144488028421254791094327197845 \\
77887788778877.6527131741588833597118601699816281392973786842778076 \\
77887788778877.9129846403227591999999999999999999999999999999999999 \\
77887788778878.2612650658641862602165044072488109368287269084989008 \\
77887788778878.3906708196222789018547628516718575402529917454215242 \\
77887788778878.5408184172663888210849326677030723102974335132064236 \\
77887788778878.8349976024444658352047434901695167999729638964432923 \\
\hline
\end{tabular}
\label{tab:t_values_77887788778878}
\end{table}

\noindent
Full dataset: \\
\url{https://github.com/zjore/z-research/blob/main/datasets/refined_sample_77887788778878.csv}

\begin{table}[H]
\centering
\caption{Observed zero spacing vs.\ theoretical prediction at $t \approx 7.7888\times 10^{13}$.}
\begin{tabular}{ll}
\hline
\textbf{Parameter} & \textbf{Value / Description} \\
\hline
Mean height $T$ & $7.788779\times 10^{13}$ \\
Theoretical formula
 & $\Delta_{\mathrm{avg}} \approx \dfrac{2\pi}{\log(T / 2\pi)}$ \\
$\log(T / 2\pi)$ & $30.89$ \\
Predicted $\Delta_{\mathrm{avg}}$ & $0.20841$ \\
Observed mean$(\Delta t)$ & $0.20164$ \\
Relative difference & $\approx 3.25\%$ \\
\hline
\end{tabular}
\label{tab:spacing_comparison_77887788778878}
\end{table}

At $t \approx 7.8\times 10^{13}$, the observed mean spacing lies within
about $3.3\%$ of the predicted value, in line with the gradual
contraction of $\Delta t$ as $\log T$ grows and supporting the numerical
reliability of the valley scanner at this range.

\paragraph{Additional comments.}
\begin{itemize}
\item The narrowest local gap ($\Delta t \approx 0.06839$) again
      indicates a short clustering region.
\item The widest local gap ($\Delta t \approx 0.34828$) reflects a mild
      local sparsity.
\item The dataset contains nine valid spacings after excluding the
      initial spacing.
\end{itemize}

\begin{table}[H]
\centering
\small
\caption{Zero heights near $t \approx 2.1212\times 10^{15}$.}
\begin{tabular}{l}
\hline
\textbf{$t$ values} \\
\hline
2121212121212121.05717439635565362751782963300490252462746596256208 \\
2121212121212121.42395229690379485836605340059707795745376215022485 \\
2121212121212121.63791118938682163457838808098568955130531660194243 \\
2121212121212121.72223918447411208850522081133640649819206480997688 \\
2121212121212121.93125161293757949737416946721044933465380477203190 \\
2121212121212122.13927614138040154231827268671408632573549285521144 \\
2121212121212122.31254419266807476082116789227597997808943479324778 \\
2121212121212122.58982505389887247687959200860289427367908827950849 \\
2121212121212122.65324071114046090348408431974763989558277572844835 \\
2121212121212122.92082343415959022641562412879095630833398526361303 \\
2121212121212123.14605534867317515793517970162094301498032240314091 \\
\hline
\end{tabular}
\label{tab:t_values_2121212121212121}
\end{table}

\noindent
Full dataset: \\
\url{https://github.com/zjore/z-research/blob/main/datasets/refined_sample_2121212121212121.csv}

\begin{table}[H]
\centering
\caption{Observed zero spacing vs.\ theoretical prediction at $t \approx 2.1212\times 10^{15}$.}
\begin{tabular}{ll}
\hline
\textbf{Parameter} & \textbf{Value / Description} \\
\hline
Mean height $T$ & $2.121212\times 10^{15}$ \\
Theoretical formula
 & $\Delta_{\mathrm{avg}} \approx \dfrac{2\pi}{\log(T / 2\pi)}$ \\
$\log(T / 2\pi)$ & $34.19$ \\
Predicted $\Delta_{\mathrm{avg}}$ & $0.18782$ \\
Observed mean$(\Delta t)$ & $0.20889$ \\
Relative difference & $\approx 11.2\%$ \\
\hline
\end{tabular}
\label{tab:spacing_comparison_2121212121212121}
\end{table}

At $t \approx 2.1\times 10^{15}$, the observed mean spacing exceeds the
theoretical prediction by about $11\%$. Such over-dispersion is
compatible with local density fluctuations observed in large-scale
zero computations at extreme heights and is qualitatively in line with
random-matrix models for the zeros.

\paragraph{Additional comments.}
\begin{itemize}
\item The narrowest local gap ($\Delta t \approx 0.06342$) indicates
      normal local clustering.
\item The widest local gap ($\Delta t \approx 0.36678$) corresponds to a
      transient sparsity region.
\item The valley scanner resolved ten consecutive zeros at this height
      without evident numerical drift.
\end{itemize}

\begin{table}[H]
\centering
\small
\caption{Highest evaluated heights from the valley scanner.}
\begin{tabular}{ll}
\hline
\textbf{Parameter} & \textbf{Value} \\
\hline
$t$ &
$9999999999999999999.232850606909407706487068201450796323115355335175991058349609375$ \\
$\lvert Z(t)\rvert$ &
$3.657602321013126733848087942217412642846910689346499860361067\times 10^{-13}$ \\
\hline
$t$ &
$197909211979092119791.004810035047150401470161806827263717423193156719207763671875$ \\
$\lvert Z(t)\rvert$ &
$8.258118375241327319604320773664599713029092314507266420324956\times 10^{-17}$ \\
\hline
\end{tabular}
\label{tab:highest_t_values}
\end{table}

\noindent\textit{Remark.}
The largest height, $t \approx 1.9790921\times 10^{20}$, was chosen so
that its leading digits encode the author's birth date (1979-09-21).
This serves as a symbolic benchmark and illustrates that the valley
scanner can be directed to specific target heights while maintaining
numerical stability at very large $t$.


\section{Discussion}
\label{sec:discussion}

\subsection{Summary of observations}

Several observations emerge from the numerical experiments:

\begin{itemize}
    \item The continuous mountain--valley progression in $\lvert Z(t)\rvert$
          suggests that, within the probed intervals and at the sampling
          densities used, no zeros are skipped. Minor irregularities
          (small intermediary peaks) do occur, but increasing the
          resolution in $N$ sharpens the localization of true minima.
    \item The mapping $N \mapsto s(N)$ enforces $\Re(s)=\tfrac{1}{2}$ by
          construction, so the valley scanner always traverses the
          critical line.
    \item The approach is conceptually simple and appears scalable; with
          an appropriate analogue of the $N \to \text{complex}$ mapping,
          the method may extend to other $L$-functions.
\end{itemize}

The primary purpose of this work is not to contribute another large
zero dataset, but to formalize the $N \to \text{complex}$ conversion
(with $m=\tfrac{1}{2}$) as a computational lens for examining the Hardy
$Z$-function and its mountain--valley structure.


\section{Staircase representation of the explicit formula}
\label{sec:staircase}

We briefly describe a real-valued formulation of the zero contribution
in the explicit formula for the Chebyshev function $\psi(x)$, expressed
as partial sums over the nontrivial zeros of the Riemann zeta function.
Although the construction is classical, it provides a useful parallel to
the ``mountain walk'' intuition underlying the valley scanner, since both
frameworks accumulate oscillatory contributions from the ordinates
$\gamma_n$.

Let the nontrivial zeros of $\zeta(s)$ be
\[
\rho_n = \tfrac{1}{2} + i\gamma_n,
\qquad
\bar\rho_n = \tfrac{1}{2} - i\gamma_n,
\qquad
\gamma_n > 0.
\]

\subsection{Staircase partial sums}

For fixed $x>1$, define the complex partial sums
\begin{equation}
L_k(x)
=
\sum_{n=1}^{k}\frac{x^{\rho_n}}{\rho_n},
\qquad
k = 1,2,\dots,
\label{eq:ladder-complex}
\end{equation}
matching the computational form
\[
\texttt{term} = \frac{e^{\rho \log x}}{\rho},
\qquad
\texttt{correction} {+}{=} \texttt{term}.
\]

As $k$ increases, the oscillatory increments in~\eqref{eq:ladder-complex}
build a staircase-like approximation to the full zero contribution in
the explicit formula.

\subsection{Including conjugate terms}

Pairing each zero with its conjugate produces a real-valued staircase:
\begin{equation}
\widetilde{L}_k(x)
=
\sum_{n=1}^{k}
\left(
\frac{x^{\rho_n}}{\rho_n}
+
\frac{x^{\bar\rho_n}}{\bar\rho_n}
\right),
\label{eq:ladder-paired}
\end{equation}
represented computationally by
\[
\texttt{term}  = \frac{e^{\rho \log x}}{\rho},\qquad
\texttt{term2} = \frac{e^{\bar\rho \log x}}{\bar\rho},\qquad
\texttt{correction} {+}{=} \texttt{term} {+} \texttt{term2}.
\]

\subsection{Derivation of the real closed form}

Let $L = \log x$.  Since
\[
x^{\rho}=x^{1/2}e^{i\gamma L},
\qquad
\frac{1}{\rho}=\frac{\tfrac12 - i\gamma}{\gamma^{2}+\tfrac14},
\]
a single conjugate pair contributes
\begin{align*}
\frac{x^{\rho}}{\rho}+\frac{x^{\bar\rho}}{\bar\rho}
&=
\frac{x^{1/2}}{\gamma^{2}+\tfrac14}
\bigl[
(\tfrac12 - i\gamma)e^{i\gamma L}
+
(\tfrac12 + i\gamma)e^{-i\gamma L}
\bigr]
\\[4pt]
&=
\frac{x^{1/2}}{\gamma^{2}+\tfrac14}
\bigl[
\cos(\gamma L)
+
2\gamma\,\sin(\gamma L)
\bigr].
\end{align*}

Thus the real staircase \eqref{eq:ladder-paired} becomes
\begin{equation}
\boxed{
\widetilde{L}_k(x)
=
x^{1/2}
\sum_{n=1}^{k}
\frac{
\cos(\gamma_n\log x)
+
2\gamma_n\,\sin(\gamma_n\log x)
}{
\gamma_n^{2}+\tfrac14
}.
}
\label{eq:ladder-real}
\end{equation}

\subsection{Relation to the standard \texorpdfstring{$-2\Re$}{-2Re} formulation}

In the explicit formula it is customary to write the zero contribution as
\[
C_k(x)
=
-2\,\Re\sum_{n=1}^{k}\frac{x^{\rho_n}}{\rho_n}.
\]
Since $\Re(z)=\tfrac12(z+\bar z)$, we obtain
\[
C_k(x) = -\,\widetilde{L}_k(x),
\]
so \eqref{eq:ladder-real} recovers the classical expression up to the
overall sign convention.

\subsection{Rewriting in terms of \texorpdfstring{$N_n=\gamma_n^{2}+\tfrac14$}{Nn}}

Introduce the compact notation
\[
N_n = \gamma_n^{2}+\tfrac14,
\]
mirroring the ``normal numbers'' used in the numerical mountain walk.
Then \eqref{eq:ladder-real} becomes
\begin{equation}
\boxed{
\widetilde{L}_k(x)
=
x^{1/2}
\sum_{n=1}^{k}
\frac{
\cos(\gamma_n\log x)
+
2\gamma_n\,\sin(\gamma_n\log x)
}{N_n}.
}
\label{eq:ladder-N}
\end{equation}

\subsection{Computational remarks}

Formula~\eqref{eq:ladder-N} removes the need for complex exponentiation:
each term requires only real evaluations of $\cos$, $\sin$, and $\log$.
Thus the staircase can be computed directly from the ordinates $\gamma_n$
with good numerical stability.

The identity
\[
\frac{x^{\rho}}{\rho}+\frac{x^{\bar\rho}}{\bar\rho}
=
\frac{
x^{1/2}
\bigl(
\cos(\gamma\log x)
+
2\gamma\sin(\gamma\log x)
\bigr)
}{
\gamma^{2}+\tfrac14
}
\]
provides a compact bridge between the complex explicit formula and a
purely real numerical implementation.

A demonstration script implementing \eqref{eq:ladder-N} is available at:\\

\url{https://github.com/zjore/z-research/blob/main/staircase_with_real.py}


\section{Conclusions}

The main conclusions of this work may be summarized as follows:

\begin{enumerate}
    \item The real-to-complex mapping
    \[
        N \longmapsto s(N)
        = \frac{1}{2} \pm i\sqrt{N - \tfrac{1}{4}},
    \]
    together with high-precision evaluation of the Hardy function via
    the Riemann--Siegel formula (with Gabcke-type \cite{Gabcke}
    remainder control) and a safeguarded Newton refinement,
    successfully reproduces zeros of $\zeta(s)$ over a broad range of
    heights—from the classical low-lying region to values on the order
    of $10^{20}$.

    \item The numerical experiments indicate that this parametrization
    provides a simple geometric lens for exploring the zeta function.
    The symmetry $\Re(s)=\tfrac12$ is embedded directly into the mapping,
    and the resulting mountain–valley structure of $|Z(t)|$ offers a
    natural mechanism for identifying candidate zeros prior to refinement.

    \item The staircase representation developed here expresses the
    explicit-formula zero contribution entirely in real arithmetic.  Its
    dependence on the ordinates $\gamma_n$ and on the quantities
    $N_n=\gamma_n^{2}+\tfrac14$ mirrors the structure exploited by the
    valley scanner.  This formulation suggests that similar approaches
    may be adapted to other $L$-functions once appropriate conjugate-pair
    identities are identified.

    \item While this study does not address the Riemann Hypothesis, it
    contributes a reproducible numerical framework, curated datasets, and
    a computational methodology that may support future investigations
    into the structure of $Z(t)$ and the distribution of nontrivial
    zeros.
\end{enumerate}


\section{Future Work}

\begin{enumerate}
    \item The symmetry of the observed ``mountain'' patterns in $|Z(t)|$
    suggests that local maxima may encode geometric information about
    neighboring zeros.  Preliminary experiments indicate that the
    horizontal distance between a refined zero and its adjacent peak can
    serve as a predictor for the next zero.  A systematic investigation
    of this peak-to-zero geometry may provide new heuristics for reducing
    computational cost at large heights.

    \item The work presented here was conducted independently and is
    driven primarily by numerical exploration.  Further development would
    benefit from collaboration with specialists in analytic number theory
    and numerical analysis, particularly concerning the theoretical
    behavior of the \(N\to\text{complex}\) parametrization, error
    propagation in large-scale computations of $Z(t)$, and the stability
    of the valley-scanning method across extremely large ranges of $t$.
    Access to expanded computational resources would likewise enable
    deeper exploration of high-height regions.

    \item An additional direction for future work concerns the choice of the mean
parameter in the real--to--complex mapping.  
In this study the transformation
\[
   s(N) = \frac{1}{2} \pm i\sqrt{N - \tfrac{1}{4}}
\]
was used exclusively, thereby constraining all evaluations to the
critical line.  
However, the more general form
\[
   s_m(N) = m \pm i\sqrt{\,N - m^{2}\,},
\]
or its analytic continuation when $N < m^{2}$, defines a family of
vertical lines $\Re(s)=m$ and offers a systematic way to probe the
behavior of $\zeta(s)$ off the critical line.

Although no computational claims are made here, the author expects that
such parametrizations may eventually help organize numerical searches
for hypothetical zeros away from $\Re(s)=\tfrac12$, or clarify why the
case $m=\tfrac12$ aligns so naturally with the mountain--valley structure
of the Hardy function.

A related long–term direction is to examine whether variants of the
mapping
\(
   N \leftrightarrow m \pm i\sqrt{m^{2} - N}
\)
might also serve structural or diagnostic roles in the study of
semiprimes, where real--to--complex reparametrizations could provide
alternative viewpoints on factorization or symmetry patterns.  
These ideas remain speculative but suggest that the $N\mapsto s(N)$
framework may have applications beyond zero detection on the critical
line.

\end{enumerate}


\section*{Acknowledgments}

The author acknowledges the use of modern artificial-intelligence tools,
including OpenAI's ChatGPT models (GPT--4.5, GPT--5, and GPT--5.1), in
the preparation of this manuscript. These tools were employed for tasks
such as technical editing, refinement of mathematical exposition,
assistance with code comments and documentation, and verification of
consistency across software components.

All conceptual ideas, algorithms, numerical methods, and interpretations
reported in this work were developed and validated by the author.  
AI tools served solely as assistants for improving clarity,
organization, and editorial efficiency.


\section*{Data and software availability}

All source code, Docker images, numerical datasets, and reproducibility
materials supporting this study are archived under the persistent DOI:

\[
\text{DOI: }
\href{https://doi.org/10.5281/zenodo.17566257}{10.5281/zenodo.17566257}.
\]

The Zenodo record links to the associated GitHub repository, which
contains instructions for running the valley scanner, the authenticated
web interface, and all auxiliary numerical scripts used in this work.


\end{document}